\documentclass[11pt]{article}
\usepackage{amsfonts}
\usepackage{mathrsfs}
\usepackage{amsfonts,amssymb}
\usepackage{amsmath,amscd}
\usepackage[all]{xy}
\usepackage{graphicx}
\usepackage{fancyhdr}
\usepackage[dvipdfm,
          colorlinks=true,
           citecolor=green,
           linkcolor=black]{hyperref}

\begin{document}
\sloppy
\newcommand{\dickebox}{{\vrule height5pt width5pt depth0pt}}
\newtheorem{Def}{Definition}[section]
\newtheorem{Bsp}{Example}[section]
\newtheorem{Prop}[Def]{Proposition}
\newtheorem{Theo}[Def]{Theorem}
\newtheorem{Lem}[Def]{Lemma}
\newtheorem{Koro}[Def]{Corollary}
\newcommand{\pd}{{\rm proj.dim\, }}
\newcommand{\id}{{\rm inj.dim\, }}
\newcommand{\fd}{{\rm fin.dim}\,}
\newcommand{\add}{{\rm add \, }}
\newcommand{\Hom}{{\rm Hom \, }}
\newcommand{\gldim}{{\rm gl.dim}\,}
\newcommand{\End}{{\rm End \, }}
\newcommand{\Ext}{{\rm Ext}}
\newcommand{\D}{\rm D \,}
\newcommand{\cpx}[1]{#1^{\bullet}}
\newcommand{\Dz}[1]{{\rm D}^+(#1)}
\newcommand{\Df}[1]{{\rm D}^-(#1)}
\newcommand{\Db}[1]{{\rm D^b}(#1)}
\newcommand{\C}[1]{{\rm C}(#1)}
\newcommand{\Cz}[1]{{\rm C}^+(#1)}
\newcommand{\Cf}[1]{{\rm C}^-(#1)}
\newcommand{\Cb}[1]{{\rm C^b}(#1)}
\newcommand{\K}[1]{{\rm K}(#1)}
\newcommand{\Kf}[1]{{\rm K}^-(#1)}
\newcommand{\Kb}[1]{{\rm K^b}(#1)}
\newcommand{\opp}{^{\rm op}}
\newcommand{\otimesL}{\otimes^{\rm\bf L}}
\newcommand{\lra}{\longrightarrow}
\newcommand{\ra}{\rightarrow}
\newcommand{\Pmodcat}[1]{#1\mbox{{\rm -Proj}}}

\newcommand{\rad}{{\rm rad \, }}
\newcommand{\Lra}{\Longrightarrow}
\newcommand{\tra}{\twoheadrightarrow}
\newcommand{\Htp}{{\rm Htp}}
\newcommand{\Mod}{{\rm Mod}}
{\Large \bf
\begin{center}  Derived equivalences for Cohen-Macaulay Auslander algebras

 \end{center}}
\medskip

\centerline{\sc Shengyong Pan}

\medskip
\abstract{Let $A$ and $B$ be Gorenstein Artin algebras of
Cohen-Macaulay finite type.  We prove that, if $A$ and $B$ are
derived equivalent, then their Cohen-Macaulay Auslander algebras are
also derived equivalent.}

\medskip {\small {\it 2000 AMS Classification}: 16G10, 18G05,
16S50; 16P10, 20C05, 18E10.

}

\section{Introduction}

Triangulated categories and derived categories were introduced by
Grothendieck and Verdier \cite{Ver}. Today, they have widely been
used in many branches: algebraic geometry, stable homotopy theory,
representation theory, etc. In the representation theory of
algebras, we will restrict our attention to the equivalences of
derived categories, that is, derived equivalences. Derived
equivalences have been shown to preserve many invariants and provide
new connection. For instance, Hochschild homology and cohomology
\cite{Ri3}, finiteness of finitistic dimension \cite{PX} have been
shown to be invariant under derived equivalences. Moreover, derived
equivalences are related to cluster categories and cluster tilting
objects \cite{BMRR}.  As is known, Rickard's Morita theory for
derived categories leaves something to be desired, though, as for
some pairs of rings, or algebras, it is currently difficult,
sometimes even impossible to verify whether there exists a tilting
complex. It is of interest to construct a new derived equivalence
from given one by finding a suitable tilting complex. Rickard
\cite{Ri2,Ri3} used tensor products and trivial extensions to get
new derived equivalences. In the recent years, Hu and Xi have
provided various techniques to construct new derived equivalences.
In \cite{HX1} they established an amazing connection between derived
equivalences and Auslander-Reiten sequences via BB-tilting modules,
and obtained derived equivalences from Auslander-Reiten triangles.
In \cite{HX3} they constructed new derived equivalences between
$\Phi$-Auslander-Yoneda algebras from a given almost $\nu$-stable
equivalence.

In \cite[Corollary 3.13]{HX3} Hu and Xi proved that, if two
representation finite self-injective Artin algebras are derived
equivalent, then their Auslander algebras are derived equivalent. In
this paper, we generalize their result and prove that, if two
Cohen-Macaulay finite Gorenstein Artin algebras are derived
equivalent, then their Cohen-Macaulay Auslander algebras are also
derived equivalent.

This paper is organized as follows. In Section 2, we review some
facts on derived categories and derived equivalences. In Section 3,
we state and prove our main result.

\section{Preliminaries}
In this section, we shall recall some definitions and notations on
derived categories and derived equivalences.

Let $\mathscr{A}$ be an abelian category. For two morphisms $\alpha:
X\ra Y$ and $\beta: Y\ra Z$, their composition is denoted by
$\alpha\beta$. An object $X\in\mathscr{A}$ is called a additive
generator for $\mathscr{A}$ if $\add(X)=\mathscr{A}$, where
$\add(X)$ is the additive subcategory of $\mathscr{A}$ consisting of
all direct summands of finite direct sums of the copies of $X$. A
complex $\cpx{X}=(X^i,d_{X}^i)$ over $\mathscr{A}$ is a sequence of
objects $X^i$ and morphisms $d_{X}^i$ in $\mathscr{A}$ of the form:
$\cdots \ra X^i\stackrel{d^i}\ra X^{i+1}\stackrel{d^{i+1}}\ra
X^{i+1}\ra\cdots$, such that $d^id^{i+1}=0$ for all
$i\in\mathbb{Z}$. If $\cpx{X}=(X^i,d_{X}^i)$ and
$\cpx{Y}=(Y^i,d_{Y}^i)$ are two complexes, then a morphism $\cpx{f}:
\cpx{X}\ra\cpx{Y}$ is a sequence of morphisms $f^i: X^i\ra Y^i$ of
$\mathscr{A}$ such that $d^i_{X}f^{i+1}=f^id^i_{Y}$ for all
$i\in\mathbb{Z}$. The map $\cpx{f}$ is called a chain map between
$\cpx{X}$ and $\cpx{Y}$. The category of complexes over
$\mathscr{A}$ with chain maps is denoted by $\C{\mathscr{A}}$. The
homotopy category of complexes over $\mathscr{A}$ is denoted by
$\K{\mathscr{A}}$ and the derived category of complexes is denoted
by $\D(\mathscr{A})$.

Let $R$ be a commutative Artin ring. And let $A$ be an Artin
$R$-algebra. We denote by $A$-mod the category of finitely generated
left $A$-modules. The full subcategory of $A$-mod consisting of
projective modules is denoted by $_A\mathcal {P}$. Recall that a
 homomorphism $f: X\ra Y$ of $A$-modules is called a radical map
 provided that for any $A$-module $Z$ and homomorphisms $g: Y\ra Z$ and $h: Z\ra
 X$, the composition $hfg$ is not an isomorphism. A complex of
 $A$-modules is called a radical complex if its differential maps
 are radical maps.
 Let $\Kb{A}$ denote the
homotopy category of bounded complexes of $A$-modules. We denote by
$\Db{A}$ by the bounded derived category of $A$-mod.

The fundamental theory on derived equivalences has been established.
Rickard \cite{Ri1} gave a Morita theory for derived categories in
the following theorem.

\begin{Theo}$\rm \cite[Therem 6.4]{Ri1}$
Let $A$ and $B$ be rings. The following conditions are equivalent.

$(i)$ $\Db{A\text{-}\Mod}$ and $\Db{B\text{-}\Mod}$ are equivalent
as triangulated categories.

$(ii)$ $\Kf{\Pmodcat{A}}$ and $\Kf{\Pmodcat{B}}$ are equivalent as
triangulated categories.

$(iii)$ $\Kb{\Pmodcat{A}}$ and $\Kb{\Pmodcat{B}}$ are equivalent as
triangulated categories.

$(iv)$ $\Kb{_{A}\mathcal {P}}$ and $\Kb{_{B}\mathcal {P}}$ are
equivalent as triangulated categories.

$(v)$ $B$ is isomorphic to $\End_{\Db{A}}(\cpx{T})$ for some complex
$\cpx{T}$ in $\Kb{_{A}\mathcal {P}}$ satisfying

         \qquad $(1)$ $\Hom_{\Db{A}}(\cpx{T},\cpx{T}[n])=0$
         for all $n\neq 0$.

         \qquad $(2)$ $\add(\cpx{T})$, the category of direct summands of
          finite direct sums of copies of $\cpx{T}$, generates
          $\Kb{_{A}\mathcal {P}}$ as a triangulated category.

Here $A$-Proj is the subcategory of $A$-Mod consisting of all
projective $A$-modules.
\end{Theo}
\noindent{\bf Remarks.} (1) The rings $A$ and $B$ are said to be
derived equivalent if $A$ and $B$ satisfy the conditions of the
above theorem. The complex $\cpx{T}$ in Theorem 2.1 is called a
tilting complex for $A$.

(2) By \cite[Corollary 8.3]{Ri1}, two Artin $R$-algebras $A$ and $B$
are said to be  derived equivalent if their derived categories
$\Db{A}$ and $\Db{B}$ are equivalent as triangulated categories. By
Theorem 2.1, Artin algebras $A$ and $B$ are derived equivalent if
and only if $B$ is isomorphic to the endomorphism algebra of a
tilting complex $\cpx{T}$. If $\cpx{T}$ is a tilting complex for
$A$, then there is an equivalence $F: \Db{A}\ra\Db{B}$ that sends
$\cpx{T}$ to $B$. On the other hand, for each derived equivalence
$F: \Db{A}\ra\Db{B}$, there is an associated tilting complex
$\cpx{T}$ for $A$ such that $F(\cpx{T})$ is isomorphic to $B$ in
$\Db{B}$.

\section{Derived equivalences for Cohen-Macaulay Auslander Algebras}
In this section, we shall prove the main result of this paper.
First, let us recall the definition of Cohen-Macaulay Auslander
algebras.
\subsection{Cohen-Macaulay Auslander algebras}

Let $A$ be an Artin algebra. Recall that $A$ is of finite
representation type provided that there are only finitely many
indecomposable finitely generated $A$-modules up to isomorphism. If
an $A$-module $X$ satisfies $\Ext_{A}^{i}(X,A)=0$ for $i>0$, then
$X$ is said to be a Cohen-Macaulay $A$-module. Denote by
$_{A}\mathcal {X}$ the category of Cohen-Macaulay $A$-modules. It is
easy to see that if $A$ is a self-injective algebra, then
$_{A}\mathcal {X}=A$-mod. By a $\Hom_{A}(-,X)$-exact sequence
$\cpx{Y}=(Y^{i},d^{i})$, we mean that the sequence $\cpx{Y}$ itself
is exact, and that $\Hom_{A}(\cpx{Y},X)$ remains to be exact. An
$A$-module $X$ is said to be Gorenstein projective if there is a
$\Hom_{A}(-,Q)$-exact sequence
$$\cdots\ra
P^{-1}\stackrel{d^{-1}}\ra P^{0}\stackrel{d^{0}}\ra
P^{1}\stackrel{d^{1}}\ra\cdots
$$ such that $X\simeq Imd^{0}$, where $P^{i}$ (for each $i$) and
$Q$ are projective $A$-modules. Denote by $A$-Gproj the subcategory
of $A$-mod consisting of Gorenstein projective $A$-modules. Note
that Gorenstein projective modules are Cohen-Macaulay $A$-modules.
Following \cite[Example 8.4(2)]{Be} an Artin algebra $A$ is said to
be of  Cohen-Macaulay finite type provided that there are only
finitely many indecomposable finitely generated Gorenstein
projective $A$-modules up to isomorphism. It is easy to see that
algebras of finite representation type are of Cohen-Macaulay finite
type. Suppose that $A$ is of Cohen-Macaulay finite type. In other
words, $A$-Gproj has an additive generator $M$, that is,
$\add(M)=A$-Gproj.

\begin{Def} $\rm\cite{Ch}$ Suppose that an Artin algebra $A$ is of
Cohen-Macaulay finite type. Let $M$ be an additive generator in
$A$-Gproj. We call $\Lambda=\End(M)$ a Cohen-Macaulay Auslander
algebra of $A$.
\end{Def}
\noindent{\bf Remark.} For a Cohen-Macaulay finite algebra $A$, its
Cohen-Macaulay Auslander algebra is unique up to Morita
equivalences.

\noindent{\bf Example.} Let $A=k[x]/(x^{2})$ and consider the Artin
algebra
$$
T_{2}(A)=\left(\begin{array}{cc}
A&A\\
0&A
\end{array}\right).
$$
Then $T_{2}(A)$ is a $1$-Gorenstein Artin algebra of Cohen-Macaulay
type \cite{FGR} or \cite{Ha2}. $T_{2}(A)$ has indecomposable
Gorenstein projective modules \cite[p.101]{BR}:
$$
M_1=\left(\begin{array}{cc}
k\\
0
\end{array}\right),
M_2=\left(\begin{array}{cc}
A\\
0
\end{array}\right),
M_3=\left(\begin{array}{cc}
A\\
A
\end{array}\right),
M_4=\left(\begin{array}{cc}
k\\
k
\end{array}\right),
M_5=\left(\begin{array}{cc}
A\\
k
\end{array}\right).
$$
Set $M=\oplus_{1\leq i\leq 5} M_{i}$. Then Cohen-Macaulay Auslander
algebra $\End_{T_{2}(A)}(M)$ of $T_{2}(A)$ is given by the following
quiver and relations $xy=0=v\alpha-yu=\alpha z=\alpha\beta\gamma$
\cite{GZ}. \vspace{-0.5cm}
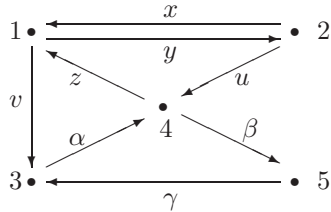
\begin{center}
\setlength{\unitlength}{1mm}
\begin{picture}  (50,20)
\put(5,10){\circle*{1}} \put(40,10){\circle*{1}}
\put(5,-10){\circle*{1}} \put(22.5,0){\circle*{1}}
\put(40,-10){\circle*{1}}

\put(2,9){1} \put(43,9){2} \put(2,-11){3} \put(22,-4){4}
\put(43,-11){5}

\put(7,9){\vector(1,0){31}} \put(5,8){\vector(0,-1){16}}
\put(38,11){\vector(-1,0){31}} \put(38,8){\vector(-2,-1){13}}
\put(7,-8){\vector(2,1){13}} \put(20,1){\vector(-2,1){13}}
\put(25,-1){\vector(2,-1){13}} \put(38,-10){\vector(-1,0){31}}

\put(22.5,6.5){$y$} \put(2,0){$v$} \put(22.5,12){$x$}
\put(32,3){$u$} \put(10,-5){$\alpha$} \put(10,3){$z$}
\put(33,-4){$\beta$} \put(22.5,-13){$\gamma$}
\end{picture}
\end{center}

\bigskip
\bigskip

\subsection{The proof of the main result}
We shall give the proof of the main result of this paper.

Suppose $A$ and $B$ are Artin algebras. Let $F: \Db{A}\lra \Db{B}$
be a derived equivalence and let $\cpx{P}$ be the tilting complex
associated to $F$. Without loss of generality, we assume that
$\cpx{P}$ is a radical complex of the following form
$$
0\ra P^{-n}\ra P^{-n+1} \ra \cdots\ra P^{-1}\ra P^{0}\ra 0.
$$
Then we have the following fact.

\begin{Lem} $\rm\cite[lemma\, 2.1]{HX1}$
Let $F: \Db{A}\lra \Db{B}$ be a derived equivalence between Artin
algebras $A$ and $B$. Then we have a tilting complex
$\bar{P}^{\bullet}$ for $B$ associated to the quasi-inverse of $F$
of the form
$$
0\ra \bar{P}^{0}\ra \bar{P}^{1} \ra \cdots\ra \bar{P}^{n-1}\ra
\bar{P}^{n}\ra 0,
$$ with the differential being radical maps.
\end{Lem}

Suppose that $\cpx{X}$ is a complex of $A$-modules. We define the
following truncations:

$\tau_{\geq 1}(\cpx{X}): \cdots\ra0\ra0\ra X^{1}\ra X^{2}\ra\cdots$,

$\tau_{\leq 0}(\cpx{X}): \cdots\ra X^{-1}\ra X^{0}\ra 0\ra0\cdots$.

Using the properties of Cohen-Macaulay $A$-modules, we can prove the
following lemma.

\begin{Lem}\label{3.3} Let $F: \Db{A}\lra \Db{B}$ be a derived equivalence between
Artin algebras $A$ and $B$, and let $G$ be the quasi-inverse of $F$.
Suppose that $\cpx{P}$ and $\bar{P}^{\bullet}$ are the tilting
complexes associated to $F$ and $G$, respectively. Then

$(i)$ For $X\in _{A}\mathcal {X}$, the complex $F(X)$ is isomorphic
in $\Db{B}$ to a radical complex $\bar{P}^{\bullet}_{X}$ of the form
$$
0\ra \bar{P}_{X}^{0}\ra \bar{P}_{X}^{1} \ra \cdots\ra
\bar{P}_{X}^{n-1}\ra \bar{P}_{X}^{n}\ra 0
$$
with $\bar{P}_{X}^{0}\in _{B}\mathcal {X}$ and $\bar{P}_{X}^{i}$
projective $B$-modules for $1\leq i\leq n$.

$(ii)$ For $Y\in_{B}\mathcal {X}$, the complex $G(Y)$ is isomorphic
in $\Db{A}$ to a radical complex $P^{\bullet}_{Y}$ of the form

$$
0\ra P_{Y}^{-n}\ra P_{Y}^{-n+1} \ra \cdots\ra P_{Y}^{-1}\ra
P_{Y}^{0}\ra 0
$$
with $P_{Y}^{-n}\in_{A}\mathcal {X}$ and $P_{Y}^{i}$ projective
$A$-modules for $-n+1\leq i\leq 0$.
\end{Lem}

\textbf{ Proof.} We just show the first case. The proof of ($ii$) is
similar to that of ($i$).

($i$) For $X\in _{A}\mathcal {X}$, by \cite[Lemma 3.1]{HX2}, we see
that the complex $F(X)$ is isomorphic in $\Db{B}$ to a complex
$\bar{P}^{\bullet}_{X}$ of the form
$$
0\ra \bar{P}_{X}^{0}\ra \bar{P}_{X}^{1} \ra \cdots\ra
\bar{P}_{X}^{n-1}\ra \bar{P}_{X}^{n}\ra 0,
$$
with $\bar{P}_{X}^{i}$ projective $B$-modules for $i>0$. We only
need to show that $\bar{P}_{X}^{0}$ is in $_{B}\mathcal {X}$. It
suffices to prove that $\End^{i}_{B}(\bar{P}_{X}^{0},B)=0$ for
$i\geq 1$. Indeed, there exists a distinguished triangle
$$\bar{P}^{+}_{X}\ra\bar{P}^{\bullet}_{X}\ra\bar{P}^{0}_{X}\ra
\bar{P}^{+}_{X}[1]$$ in $\Kb{B}$, where $\bar{P}^{+}_{X}$ denotes
the complex $\tau_{\geq 1}(\cpx{\bar{P}_{X}})$. For each
$i\in\mathbb{Z}$, applying the functor $\Hom_{\Db{B}}(-,B[i])$ to
the above distinguished triangle, we get an exact sequence
\begin{eqnarray*}
\cdots\ra\Hom_{\Db{B}}(\bar{P}^{+}_{X}[1],B[i])\ra
\Hom_{\Db{B}}(\bar{P}^{0}_{X},B[i])\ra
\Hom_{\Db{B}}(\bar{P}^{\bullet}_{X},B[i])\\\ra
\Hom_{\Db{B}}(\bar{P}^{+}_{X},B[i])\ra\cdots. \end{eqnarray*} On the
other hand, $\Hom_{\Db{B}}(\bar{P}^{+}_{X},B[i])\simeq
\Hom_{\Kb{B}}(\bar{P}^{+}_{X},B[i])=0$ for $i\geq0$. By \cite[lemma
2.1]{PX} and $\End^{i}_{A}(X,A)=0$ for $i\geq 1$, we get
$\Hom_{\Db{B}}(\bar{P}^{\bullet}_{X},B[i])\simeq\Hom_{\Db{A}}(X,P^{\bullet}[i])=0$
for all $i\geq 1$. Consequently, we get
$\Hom_{\Db{B}}(\bar{P}^{0}_{X},B[i])=0$ for all $i\geq 1$ by the
above exact sequence. Therefore,
$$
\End^{i}_{B}(\bar{P}_{X}^{0},B)\simeq
\Hom_{\Db{B}}(\bar{P}^{0}_{X},B[i])=0, \;\;\;\text{for}\;\;\; i\geq
1.$$ This implies that $\bar{P}_{X}^{0}\in_{B}\mathcal {X}$.
$\square$

Now we give a lemma, which is useful in the following argument.

\begin{Lem} Let $A$ be an Artin algebra and $f: X\ra Y$ a
homomorphism of $A$-modules with $X,Y\in_{A}\mathcal {X}$. Suppose
$\cpx{Q}$ is a complex in $\Kb{_{A}\mathcal {P}}$. If $f$ factors
through $\cpx{Q}$ in $\Db{A}$, then $f$ factors through a projective
$A$-module.
\end{Lem}

\textbf{Proof.} There is a distinguished triangle
$$\tau_{\leq
0}(\cpx{Q})\ra\tau_{\geq
1}(\cpx{Q})\stackrel{a}\ra\cpx{Q}\stackrel{b}\ra\tau_{\leq
0}(\cpx{Q})[1]\quad\text{in}\quad \Db{A}.$$ Suppose that $f=gh$,
where $g: X\ra \cpx{Q}$ and $h:\cpx{Q}\ra Y$. Since
$\Hom_{\Db{A}}(\tau_{\geq
1}(\cpx{Q}),Y)\simeq\Hom_{\Kb{A}}(\tau_{\geq 1}(\cpx{Q}),Y)=0$, it
follows that $ah=0$. Then there is a map $x: \tau_{\leq
0}(\cpx{Q})[1]\ra Y$, such that $h=bx$. Thus, we get $f=gbx$. Now,
it is sufficient to show that $f$ factors through $\tau_{\leq
0}(\cpx{Q})$. Consider the following distinguished triangle
$$Q^{0}\stackrel{c}\ra\tau_{\leq 0}(\cpx{Q})\stackrel{d}\ra\tau_{\leq
-1}(\cpx{Q})\ra Q^{0}[1]\quad\text{in}\quad\Db{A}.$$ Note that
$\Ext_{A}^{i}(X,A)=0$ for $i\geq 1$. By \cite[Lemma 2.1]{PX}, we
have $\Hom_{\Db{A}}(X,\tau_{\leq -1}(\cpx{Q}))=0$. Thus, we get
$gbd=0$. Then there is a morphism $u: X\ra Q^{0}$ such that $gb=uc$.
Consequently, $f=ucx$, which implies that $f$ factors through a
projective $A$-module $Q^{0}$. $\square$

Choose an $A$-module $X\in_{A}\mathcal {X}$, by Lemma 3.3, we know
that $F(X)$ is isomorphic to a radical complex of the form
$$
0\ra \bar{P}_{X}^{0}\ra \bar{P}_{X}^{1} \ra \cdots\ra
\bar{P}_{X}^{n-1}\ra \bar{P}_{X}^{n}\ra 0
$$
such that $\bar{P}_{X}^{0}\in _{B}\mathcal {X}$ and
$\bar{P}_{X}^{i}$ are projective $B$-modules for $1\leq i\leq n$. In
the following, we try to define a functor $\underline{F}:
\underline{_{A}\mathcal {X}}\ra \underline{_{B}\mathcal {X}}$.
\begin{Prop}
Let $F:\Db{A}\lra \Db{B}$ be a derived equivalence. Then there is an
additive functor $\underline{F}: \underline{_{A}\mathcal
{X}}\ra\underline{_{B}\mathcal {X}}$ sending $X$ to
$\bar{P}_{X}^{0}$, such that the following diagram
$$\xymatrix{
      \underline{_{A}\mathcal {X}}\ar[r]^(.35){{\rm can}}\ar[d]_{\underline{F}} &
      \Db{A}/\Kb{_{A}\mathcal {P}}\ar[d]^{{F}}\\
      \underline{_{B}\mathcal {X}}\ar[r]^(.35){{\rm can}} & \Db{B}/\Kb{_{B}\mathcal {P}}
    }$$
   is commutative up to natural isomorphism.
\end{Prop}

\textbf{ Proof.} The idea of the proof is similar to that of
\cite[Proposition 3.4]{HX1}. For convenience, we give the details
here.

For each $f:X\ra Y$ in $_{A}\mathcal {X}$, we denote by
$\underline{f}$ the image of $f$ in $\underline{_{A}\mathcal {X}}$.
By Lemma \ref{3.3}, we have a distinguished triangle
$$\bar{P}^{+}_{X}\stackrel{i_X}\ra
F(X)\stackrel{j_X}\ra\bar{P}^{0}_{X}\stackrel{m_X}\ra\bar{P}^{+}_{X}[1]\quad
\text{in}\quad \Db{B}.$$ Moreover, for each $f:X\ra Y$ in
$_{A}\mathcal {X}$, there is a commutative diagram
$$\xymatrix{
\bar{P}^{+}_{X}\ar^{i_{X}}[r]\ar^{\alpha_{f}}[d]
&F(X)\ar^{j_{X}}[r]\ar^{F(f)}[d]
  & \bar{P}^{0}_{X}\ar^{m_{X}}[r]\ar^{\beta_{f}}[d]&\bar{P}^{+}_{X}[1]
  \ar^{\alpha_{f[1]}}[d]\\
\bar{P}^{+}_{Y}\ar^{i_{Y}}[r] & F(Y)\ar^{j_{Y}}[r]
  & \bar{P}^{0}_{Y}\ar^{m_{Y}}[r]& \bar{P}^{+}_{Y}[1] .}$$
Since $\Hom_{\Db{B}}(\bar{P}^{+}_{X},\bar{P}^{0}_{Y})\simeq
\Hom_{\Kb{B}}(\bar{P}^{+}_{X},\bar{P}^{0}_{Y})=0$, it follows that
$i_{X}F(f)j_{Y}=0$. Then there exists a homomorphism $\alpha_{f}:
\bar{P}^{+}_{X}\ra \bar{P}^{+}_{Y}$. Note that $B$-mod is fully
embedding into $\Db{B}$, hence $\beta_{f}$ is a morphism of
$B$-modules. If there is another morphism $\beta'_{f}$ such that
$j_{X}\beta'_{f}=F(f)j_{Y}$, then $j_{X}(\beta_{f}-\beta'_{f})=0$.
Thus $\beta_{f}-\beta'_{f}$ factors through $\bar{P}^{+}_{X}[1]$,
which implies that $\beta_{f}-\beta'_{f}$ factors through a
projective $B$-module by Lemma 3.4. Therefore, the morphism
$\bar{\beta_{f}}$ is uniquely determined by $f$.

Let $f: X\ra Y$ and $g: Y\ra Z$ be morphisms in $_{A}\mathcal {X}$.
Then we have $F(fg)j_{Z}=j_{X}\beta_{fg}$ and
$F(fg)j_{Z}=j_{X}\beta_{f}\beta_{g}$. By the uniqueness of
$\underline{\beta_{fg}}$, we have
$\underline{\beta_{fg}}=\underline{\beta_{f}}$
$\underline{\beta_{g}}$. Moreover, if $f$ factors through a
projective $A$-module, then $\beta_{f}$ also factors through a
projective $B$-module.

For each $X\in_{A}\mathcal {X}$, we define
$\underline{F}(X):=\bar{P}_{X}^{0}$. Set
$\underline{F}(\underline{f})=\underline{\beta_{f}}$, for each
$\underline{f}\in\Hom_{\underline{\mathcal {X}_{A}}}(X,Y)$. Then $F$
is well-defined and an additive functor.

To complete the proof of the lemma, it remains to show that $j_X:
F(X)\ra \underline{F}(X)$ is a natural isomorphism in
$\Db{B}/\Kb{_{B}\mathcal {P}}$. Since $\bar{P}^{+}_{X}$ is in
$\Kb{_{B}\mathcal {P}}$, then $j_X: F(X)\stackrel{\simeq}\ra
\underline{F}(X)$ in $\Db{B}/\Kb{_{B}\mathcal {P}}$. The following
commutative diagram
$$\xymatrix{
F(X)\ar^{j_{X}}[r]\ar^{F(f)}[d]
  & \underline{F}(X)=\bar{P}^{0}_{X}\ar^{\beta_{f}}[d]\\
F(Y)\ar^{j_{Y}}[r]
  & \underline{F}(Y)=\bar{P}^{0}_{Y} }$$
shows that $j_X: F(X)\ra \underline{F}(X)$ is a natural isomorphism
in $\Db{B}/\Kb{_{B}\mathcal {P}}$. $\square$

The following lemma is quoted from \cite{HX1} which will be used
frequently.

\begin{Lem} $\rm \cite[Lemma\, 2.2]{HX1}$
Let $R$ be an arbitrary ring, and let $R$-Mod be the category of all
left $A$-modules. Suppose $X^{\bullet}$ is a bounded above complex
and $Y^{\bullet}$ is a bounded below complex over $R$-Mod. Let $m$
be an integer. If $X^{i}$ is projective for all $i>m$ and $Y^{j}=0$
for all $j<m$, then
$\Hom_{\K{R\text{-}\Mod}}(X^{\bullet},Y^{\bullet})\simeq\Hom_{\D(R\text{-}Mod)}(X^{\bullet},Y^{\bullet})$.
\end{Lem}

Let $A$ be an Artin algebra and let $X$ be in $_{A}\mathcal {X}$
which is not a projective $A$-module. Set $\Lambda=\End_{A}(A\oplus
X)$, $N=B\oplus\underline{F}(X)$ and $\Gamma=\End_{B}(N)$. Let
$\bar{T}^{\bullet}$ be the complex
$\bar{P}^{\bullet}\oplus\bar{P}^{\bullet}_{X}$. Then
$\bar{T}^{\bullet}$ is in $\Kb{\add_{B}N}$.

The proof of the following lemma is different from \cite[Lemma
3.6]{HX3}, and in fact extends Hu and Xi's original methods for the
self-injective case.

\begin{Lem} Keep the notations above. We have the following statements.

$(1)$
$\Hom_{\Kb{\add_{B}N}}(\bar{T}^{\bullet},\bar{T}^{\bullet}[i])=0$
for $i\neq0$.

$(2)$ $\add\bar{T}^{\bullet}$ generates $\Kb{\add_{B}N}$ as a
triangulated category.
\end{Lem}

\textbf{Proof.} (1) Decompose the complex $\bar{T}^{\bullet}$ as
$\bar{P}^{\bullet}\oplus\bar{P}^{\bullet}_{X}$. Then we have the
following isomorphisms
\begin{eqnarray*}
\Hom_{\Kb{B}}(\bar{T}^{\bullet},\bar{T}^{\bullet}[i])\simeq
\Hom_{\Kb{B}}(\bar{P}^{\bullet}\oplus
\bar{P}^{\bullet}_{X},(\bar{P}^{\bullet}\oplus\bar{P}^{\bullet}_{X})[i])
\simeq\Hom_{\Kb{B}}(\bar{P}^{\bullet},\bar{P}^{\bullet}[i])\oplus\\
\Hom_{\Kb{B}}(\bar{P}^{\bullet},\bar{P}^{\bullet}_{X}[i])
\oplus\Hom_{\Kb{B}}(\bar{P}^{\bullet}_{X},\bar{P}^{\bullet}[i])\oplus
\Hom_{\Kb{B}}(\bar{P}^{\bullet}_{X},
\bar{P}^{\bullet}_{X}[i]).\end{eqnarray*}

The proof falls naturally into three parts.

(a) Since $\bar{P}^{\bullet}$ is a tilting complex over $B$, we have
$\Hom_{\Kb{B}}(\bar{P}^{\bullet},\bar{P}^{\bullet}[i])=0$ for all
$i\neq0$. Furthermore,
$$\Hom_{\Kb{B}}(\bar{P}^{\bullet},\bar{P}^{\bullet}_{X}[i])\simeq
\Hom_{\Db{B}}(\bar{P}^{\bullet},\bar{P}^{\bullet}_{X}[i])\simeq
\Hom_{\Db{A}}(A,X[i])=0\quad \text{for all}\quad i\neq0.$$

(b) We claim that
$\Hom_{\Kb{B}}(\bar{P}^{\bullet}_{X},\bar{P}^{\bullet}[i])=0$ for
$i\neq0$.

Indeed, applying the functors $\Hom_{\K{B}}(-,\bar{P}^{\bullet}[i])$
and $\Hom_{\Db{B}}(-,\bar{P}^{\bullet}[i])$ to the distinguished
triangle $
\bar{P}^{+}_{X}\ra\bar{P}^{\bullet}_{X}\ra\bar{P}^{0}_{X}\ra
\bar{P}^{+}_{X}[1]$ in $\Kb{B}$, we obtain the following commutative
diagram
$$\xymatrix{
\Hom_{\Kb{B}}(\bar{P}^{+}_{X}[1],\bar{P}^{\bullet}[i])\ar[r]\ar^{\simeq}[d]&
\Hom_{\Kb{B}}(\bar{P}^{0}_{X},\bar{P}^{\bullet}[i])\ar[r]\ar[d]
&\Hom_{\Kb{A}}(\bar{P}^{\bullet}_{X},\bar{P}^{\bullet}[i])\ar[r]\ar[d]
  & \cdots \\
  \Hom_{\Db{B}}(\bar{P}^{+}_{X}[1],\bar{P}^{\bullet}[i])\ar[r]&
\Hom_{\Db{A}}(\bar{P}^{0}_{X},\bar{P}^{\bullet}[i])\ar[r] &
\Hom_{\Db{A}}(\bar{P}^{\bullet}_{X},\bar{P}^{\bullet}[i])\ar[r]
  & \cdots.}$$
Note that
$$
\Hom_{\Kb{B}}(\bar{P}^{0}_{X},\bar{P}^{\bullet})\simeq
\Hom_{\Db{B}}(\bar{P}^{0}_{X},\bar{P}^{\bullet}).
$$
Indeed, since $\cpx{P}$ is a bounded complex and Lemma 3.6, it
suffices to show that for the complex $\bar{P}^{\bullet}$ of length
$2$ of the form $0\ra \bar{P}^{0}\ra \bar{P}^{1}\ra0$, we get
$$
\Hom_{\Kb{B}}(\bar{P}^{0}_{X},\bar{P}^{\bullet})\simeq
\Hom_{\Db{B}}(\bar{P}^{0}_{X},\bar{P}^{\bullet}).
$$
In this case, we have a distinguished triangle
$$\bar{P}^{1}[-1]\ra\bar{P}^{\bullet}\ra\bar{P}^{0} \ra
\bar{P}^{1}\quad \text{in}\quad \Kb{B}.$$ Applying the functors
$\Hom_{\Kb{B}}(\bar{P}^{0}_{X},-)$ and
$\Hom_{\Db{B}}(\bar{P}^{0}_{X},-)$ to the distinguished triangle
$\bar{P}^{1}[-1]\ra\bar{P}^{\bullet}\ra\bar{P}^{0} \ra \bar{P}^{1}$,
we obtain the following commutative diagram
$$
\xymatrix{
\Hom_{\Kb{B}}(\bar{P}^{0}_{X},\bar{P}^{1}[-1])\ar[r]\ar^{\simeq}[d]
&\Hom_{\Kb{B}}(\bar{P}^{0}_{X},\bar{P}^{\bullet})\ar[r]\ar[d]
  & \Hom_{\Kb{B}}(\bar{P}^{0}_{X},\bar{P}^{0})\ar[r]\ar^{\simeq}[d]&
  \Hom_{\Kb{B}}(\bar{P}^{0}_{X},\bar{P}^{1}) \ar^{\simeq}[d]\\
 \Hom_{\Db{B}}(\bar{P}^{0}_{X},\bar{P}^{1}[-1])\ar[r] &
\Hom_{\Db{B}}(\bar{P}^{0}_{X},\bar{P}^{\bullet})\ar[r]
  & \Hom_{\Db{B}}(\bar{P}^{0}_{X},\bar{P}^{0})\ar[r]&
  \Hom_{\Db{B}}(\bar{P}^{0}_{X},\bar{P}^{1}). }
$$
Since $\Hom_{\Db{B}}(\bar{P}^{0}_{X},\bar{P}^{1}[-1])=0$ and
$\Hom_{\Kb{B}}(\bar{P}^{0}_{X},\bar{P}^{1}[1])=0$, we conclude that
that, for $i\neq0$, we have
$$\Hom_{\Kb{B}}(\bar{P}^{\bullet}_{X},\bar{P}^{\bullet}[i])\simeq
\Hom_{\Db{B}}(\bar{P}^{\bullet}_{X},\bar{P}^{\bullet}[i])\simeq
\Hom_{\Db{A}}(X,A[i])=0.$$

(c) We claim that $\Hom_{\Kb{B}}(\bar{P}^{\bullet}_{X},
\bar{P}^{\bullet}_{X}[i])=0$ for $i\neq0$.

Indeed, it follows that $\Hom_{\Kb{B}}(\bar{P}^{\bullet}_{X},
\bar{P}^{\bullet}_{X}[i])=0$ for $i<0$ by Lemma 3.6. It suffices to
show that
$\Hom_{\Kb{B}}(\bar{P}^{\bullet}_{X},\bar{P}^{\bullet}_{X}[i])=0$
for $i>0$. Note that there is a distinguished triangle $$
(\star)\quad\quad\bar{P}^{+}_{X}\ra\bar{P}^{\bullet}_{X}\ra\bar{P}^{0}_{X}\ra
\bar{P}^{+}_{X}[1]\quad \text{in}\quad \Kb{B},
\text{\;where\;}\bar{P}^{+}_{X}\text{\;denotes the complex\;}
\tau_{\geq 1}(\cpx{\bar{P}_{X}}) .$$ Applying the functor
$\Hom_{\Kb{B}}(\bar{P}^{\bullet}_{X},-)$ to ($\star$), we get a long
exact sequence
$$
\cdots\ra\Hom_{\Kb{B}}(\bar{P}^{\bullet}_{X},\bar{P}^{+}_{X}[i])\ra
\Hom_{\Kb{B}}(\bar{P}^{\bullet}_{X},\bar{P}^{\bullet}_{X}[i])\ra
\Hom_{\Kb{B}}(\bar{P}^{\bullet}_{X},\bar{P}^{0}_{X}[i])\ra\cdots
(\star\star).
$$
From the distinguished triangle
$\bar{P}^{+}_{X}\ra\bar{P}^{\bullet}_{X}\ra\bar{P}^{0}_{X}\ra
\bar{P}^{+}_{X}[1]$, we conclude that $ H^{i}(G(\bar{P}^{+}_{X}))=0$
for $i>1$ and $G(\bar{P}^{+}_{X})$ is a radical complex
$Q^{\bullet}_{X}$ of the form
$$
\cdots\ra Q^{-1}_{X}\ra Q^{0}_{X}\ra Q^{1}_{X}\ra 0.
$$
Applying the functors $\Hom_{\Kb{B}}(-,\bar{P}^{+}_{X}[i])$ and
$\Hom_{\Kb{B}}(-,\bar{P}^{+}_{X}[i])$ to ($\star$) again, we have
the following commutative diagram
$$\xymatrix{
\Hom_{\Kb{B}}(\bar{P}^{+}_{X}[1],\bar{P}_{X}^{+}[i])\ar[r]\ar^{\simeq}[d]&
\Hom_{\Kb{B}}(\bar{P}^{0}_{X},\bar{P}_{X}^{+}[i])\ar[r]\ar^{\simeq}[d]
&\Hom_{\Kb{B}}(\bar{P}^{\bullet}_{X},\bar{P}_{X}^{+}[i])\ar[r]\ar[d]
  & \cdots \\
  \Hom_{\Db{B}}(\bar{P}^{+}_{X}[1],\bar{P}_{X}^{+}[i])\ar[r]&
\Hom_{\Db{B}}(\bar{P}^{0}_{X},\bar{P}_{X}^{+}[i])\ar[r] &
\Hom_{\Db{B}}(\bar{P}^{\bullet}_{X},\bar{P}_{X}^{+}[i])\ar[r]
  & \cdots.}$$
Therefore,
\begin{eqnarray*}
\Hom_{\Kb{B}}(\bar{P}^{\bullet}_{X},\bar{P}^{+}_{X}[i])\simeq
\Hom_{\Db{B}}(\bar{P}^{\bullet}_{X},\bar{P}^{+}_{X}[i])\simeq
\Hom_{\Db{A}}(G(\bar{P}^{\bullet}_{X}),G(\bar{P}^{+}_{X}[i]))\\\simeq
\Hom_{\Db{A}}(X,G(\bar{P}^{+}_{X})[i]).
\end{eqnarray*}
By \cite[lemma 2.1]{PX}, it follows that
$\Hom_{\Db{A}}(X,G(\bar{P}^{+}_{X})[i])=0$ for all $i>1$.
Consequently,
$\Hom_{\Kb{B}}(\bar{P}^{\bullet}_{X},\bar{P}^{+}_{X}[i])=0$ for
$i>1$. Since
$\Hom_{\Kb{B}}(\bar{P}^{\bullet}_{X},\bar{P}^{0}_{X}[i])=0$ for
$i>0$ by shifting, it follows that
$\Hom_{\Kb{B}}(\bar{P}^{\bullet}_{X},\bar{P}^{\bullet}_{X}[i])=0$
for $i>1$ by the long exact sequence $(\star\star)$. It remains to
prove that
$\Hom_{\Kb{B}}(\bar{P}^{\bullet}_{X},\bar{P}^{\bullet}_{X}[1])=0$.
To get $\Hom_{\Kb{B}}(\bar{P}^{\bullet}_{X},
\bar{P}^{\bullet}_{X}[1])=0$, it suffices to show that the map
$$
(\maltese)\quad\quad
\Hom_{\Kb{B}}(\bar{P}^{\bullet}_{X},\bar{P}^{0}_{X})\ra
\Hom_{\Kb{B}}(\bar{P}^{\bullet}_{X},\bar{P}^{+}_{X}[1])\quad
\text{is \; surjective}.
$$

From the above argument, we have the following commutative diagram
$$
\xymatrix{ Q^{\bullet}_{X} \ar^{a}[r]\ar_{\simeq}[d]  &
X\ar[r]\ar_{\simeq}[d]
& M(a) \ar[r]\ar_{\simeq}[d]& Q^{\bullet}_{X}[1]  \ar_{\simeq}[d]  \\
 G(\bar{P}^{+}_{X}) \ar[r]  &  G(\bar{P}^{\bullet}_{X})\ar[r]
  &  G(\bar{P}^{0}_{X})\ar[r]&
 G(\bar{P}^{+}_{X})[1]}
$$in $\Db{A}$, where all the vertical maps are isomorphisms, and
the morphism $a$ is chosen in $\Kb{A}$ such that the first square is
commutative. Applying the functor $\Hom_{\Kb{A}}(X,-)$ to the first
horizontal distinguished triangle, we get an exact sequence
$$\Hom_{\Kb{A}}(X,M(a))\ra \Hom_{\Kb{A}}(X,Q^{\bullet}_{X}[1])\ra
0,\quad \text{since}\quad \Hom_{\Kb{A}}(X,X[1])=0.$$ We have the
following formulas.
\begin{eqnarray*}
\Hom_{\Kb{B}}(\bar{P}^{\bullet}_{X},\bar{P}^{0}_{X})\stackrel{(\ast)}\simeq
\Hom_{\Db{B}}(\bar{P}^{\bullet}_{X},\bar{P}^{0}_{X})\simeq
\Hom_{\Db{A}}(G(\bar{P}^{\bullet}_{X}),G(\bar{P}^{0}_{X}))\\\simeq
\Hom_{\Db{A}}(X,M(a))
\end{eqnarray*}
and
\begin{eqnarray*}
\Hom_{\Kb{B}}(\bar{P}^{\bullet}_{X},\bar{P}^{+}_{X}[1])\stackrel{(\ast\ast)}\simeq
\Hom_{\Db{B}}(\bar{P}^{\bullet}_{X},\bar{P}^{+}_{X}[1])\simeq
\Hom_{\Db{A}}(G(\bar{P}^{\bullet}_{X}),G(\bar{P}^{+}_{X})[1])\\\simeq
\Hom_{\Db{A}}(X,Q^{\bullet}_{X}[1]). \end{eqnarray*} The
isomorphisms $(\ast)$ and $(\ast\ast)$ are deduced by Lemma 3.6.
Then we have the following commutative diagram
$$
\xymatrix{
\Hom_{\Kb{B}}(\bar{P}^{\bullet}_{X},\bar{P}_{X}^{0})\ar[r]\ar^{\simeq}[d]&
\Hom_{\Kb{B}}(\bar{P}^{\bullet}_{X},\bar{P}_{X}^{+})\ar^{\simeq}[d] \\
  \Hom_{\Db{A}}(X,M(a))\ar[r]&
\Hom_{\Db{A}}(X,Q^{+}[1]) .}
$$
From the above diagram, to show the map ($\maltese$) is surjective,
it is sufficient to show the map
$$
\Hom_{\Db{A}}(X,M(a))\ra \Hom_{\Db{A}}(X,\cpx{Q}_{X}[1])\quad
\text{is \; surjective}.
$$
Applying the functor $\Hom(X,-)$ and $\Hom(X,-)$ to the
distinguished triangle $\cpx{Q}\ra X\ra M(a)\ra \cpx{Q}[1]$, we get
the following commutative diagram
$$\xymatrix{
\Hom_{\Kb{A}}(X,M(a))\ar[r]\ar[d]
  & \Hom_{\Kb{A}}(X,Q_{X}^{\bullet}[1])\ar[r]\ar[d]& 0 \ar[d]\\
\Hom_{\Db{A}}(X,M(a))\ar[r]
  & \Hom_{\Db{A}}(X,Q_{X}^{\bullet}[1])\ar[r]& \Hom_{\Db{A}}(X,X[1]).}$$
Thus, to get the map
$$
\Hom_{\Db{A}}(X,M(a))\ra \Hom_{\Db{A}}(X,\cpx{Q}_{X}[1])\quad
\text{is\; surjective},
$$it suffices to show the following isomorphisms
$$(i)\quad\Hom_{\Db{A}}(X,Q^{\bullet}_{X}[1])\simeq\Hom_{\Kb{A}}(X,Q^{\bullet}_{X}[1])$$
and
$$
(ii)\quad\Hom_{\Db{A}}(X,M(a))\simeq\Hom_{\Kb{A}}(X,M(a)).
$$
Firstly, we show that
$$
(i)\quad\Hom_{\Db{A}}(X,Q^{\bullet}_{X}[1])\simeq\Hom_{\Kb{A}}(X,Q^{\bullet}_{X}[1]).$$
Indeed, it suffices to show that for the complex $Q^{\bullet}_{X}$
of the form $ 0\ra Q_{X}^{-1}\ra Q_{X}^{0}\ra0$, we get ($i$). There
is a distinguished triangle $$(\clubsuit)\quad\quad Q_{X}^{-1}\ra
Q_{X}^{0}\ra Q_{X}^{\bullet}\ra Q_{X}^{-1}[1]\quad \text{in} \quad
\Kb{A}.$$ Applying the functors $\Hom_{\Kb{A}}(X,-)$,
$\Hom_{\Db{A}}(X,-)$ to ($\clubsuit$), we obtain the following
commutative diagram
$$\xymatrix{
\Hom_{\Kb{A}}(X,Q_{X}^{-1})\ar[r]\ar^{\simeq}[d] &
\Hom_{\Kb{A}}(X,Q_{X}^{0})\ar[r]\ar^{\simeq}[d]
  & \Hom_{\Kb{A}}(X,Q_{X}^{\bullet})\ar[r]\ar[d]& \Hom_{\Kb{A}}(X,P^{-1}[1]) \ar^{\simeq}[d]\\
 \Hom_{\Db{A}}(X,Q_{X}^{-1})\ar[r] &
\Hom_{\Db{A}}(X,Q_{X}^{0})\ar[r]
  & \Hom_{\Db{A}}(X,Q_{X}^{\bullet})\ar[r]& \Hom_{\Db{A}}(X,P^{-1}[1]).}$$ Since
$\End^{i}_{A}(X,A)=0$ for $i\geq 1$, it follows that
$\Hom_{\Db{A}}(X,P^{-1}[1])=0$. Moreover,\\
$\Hom_{\Kb{A}}(X,P^{-1}[1])=0$. We thus get
$\Hom_{\Db{A}}(X,Q^{\bullet}_{X}[1])\simeq\Hom_{\Kb{A}}(X,Q^{\bullet}_{X}[1])$.
Next, we prove that
$$
(ii)\quad\Hom_{\Db{A}}(X,M(a))\simeq\Hom_{\Kb{A}}(X,M(a)).
$$
Indeed, there exists a distinguished triangle
$$
(\spadesuit)\quad\quad M(a)^{0}\ra M(a)\ra M(a)^{-}\ra
M(a)^{0}[1]\quad \text{in}\quad \Kb{A},$$ where $M(a)^{-}$ denotes
the truncated complex $\tau_{\leq-1}(M(a))$. Applying the
homological functors $\Hom_{\Kb{A}}(X,-)$ and $\Hom_{\Db{A}}(X,-)$
to ($\spadesuit$), we obtain the following commutative diagram
$$\xymatrix{
_{\Kb{A}}(X,M(a)^{-}[-1])\ar[r]\ar^{\simeq}[d]&
_{\Kb{A}}(X,M(a)^{0})\ar[r]\ar^{\simeq}[d] &
_{\Kb{A}}(X,M(a))\ar[r]\ar[d]
  & _{\Kb{A}}(X,M(a)^{-}) \ar[d]\\
_{\Db{A}}(X,M(a)^{-}[-1])\ar[r]& _{\Db{A}}(X,M(a)^{0})\ar[r] &
_{\Db{A}}(X,M(a))\ar[r]
  & _{\Db{A}}(X,M(a)^{-}).}$$
According to ($i$), we have
$\Hom_{\Kb{A}}(X,M(a)^{-})\simeq\Hom_{\Db{A}}(X,M(a)^{-})=0$.
Therefore, $\Hom_{\Db{A}}(X,M(a))\simeq\Hom_{\Kb{A}}(X,M(a))$.

From the above argument, we have shown that
$\Hom_{\Kb{B}}(\bar{P}^{\bullet}_{X}, \bar{P}^{\bullet}_{X}[i])=0$
for $i\neq0$. Since $\Kb{B}$ is a full subcategory of
$\Kb{\add_{B}N}$, it follows that
$\Hom_{\Kb{\add_{B}N}}(\bar{P}^{\bullet}_{X},
\bar{P}^{\bullet}_{X}[i])=0$ for $i\neq0$.

(2) Since $\bar{P}^{\bullet}$ is a tilting complex for $B$, we see
that $\add\bar{P}^{\bullet}$ generates $\Kb{\add_{B}B}$ as
triangulated category. All the terms of $\bar{P}^{+}_{X}$ are in
$\add_{B}B$. From the distinguished triangle
$$\bar{P}^{+}_{X}\ra\bar{P}^{\bullet}_{X}\ra\bar{P}^{0}_{X}\ra
\bar{P}^{+}_{X}[1],$$ it follows that $\bar{P}^{0}_{X}$ is in the
triangulated subcategory generated by
$\add(\bar{P}^{\bullet}\oplus\bar{P}^{\bullet}_{X})$. Therefore,
$\add\bar{T}^{\bullet}$ generates $\Kb{\add_{B}N}$ as a triangulated
category. $\square$

\begin{Prop}
The complex $\Hom(N,\bar{T}^{\bullet})$ is a tiling complex over
$\Gamma$ with the endomorphism
$\End(\Hom(N,\bar{T}^{\bullet}))\simeq\Lambda$. In particular, Artin
algebras $\Lambda$ and $\Gamma$ are derived equivalent associated
with the tilting complex $\Hom(N,\bar{T}^{\bullet})$.
\end{Prop}

\textbf{Proof.} We have an equivalence of categories
$$
\Hom_{B}(N,-): \add_{B} N\stackrel{\simeq}\lra _{\Gamma}\mathcal
{P}.
$$
We thus get an equivalence of triangulated categories induced by
$\Hom_{B}(N,-)$ as follows
$$
\Kb{\add_{B} N}\stackrel{\simeq}\lra \Kb{_{\Gamma}\mathcal {P}}.
$$
Then $\Hom(N,\bar{T}^{\bullet})\in \Kb{_{\Gamma}\mathcal {P}}$. By
the Lemma 3.7, we see that $\add\Hom(N,\bar{T}^{\bullet})$ generates
$\Kb{_{\Gamma}\mathcal {P}}$ as a triangulated category, and $
\End(\Hom(N,\bar{T}^{\bullet}))\simeq\End(\bar{T}^{\bullet})\simeq\Lambda
$. $\square$

We have the following lemma, its proof is due to Happel \cite[Lemma
4.4]{Ha2}.

\begin{Lem} Suppose that $\id_{A}A<\infty$. Then the following statements are
equivalent.

$(i)$ $\pd_{A}(\D(A_{A}))<\infty$.

$(ii)$ For $X\in _{A}\mathcal {X}$, there exists an exact sequence $
0\ra X\ra P\ra X^{'}\ra 0$, with $X^{'}\in_{A}\mathcal {X}$ and $P$
a projective $A$-module.

$(iii)$ If $X\in_{A}\mathcal {X}$ satisfies $\id_{A}X<\infty$, then
$X$ is a projective $A$-module.
\end{Lem}
Recall that an Artin algebra $A$ is called Gorenstein if the regular
module A has finite injective dimension on both sides. If $A$ is a
Gorenstein algebra, then it follows from Lemma $7.2.8$ that
$_{A}\mathcal {X}=A$-Gproj, and that $_{A}\mathcal {X}$ is a
Frobenius category and its category $\underline{_{A}\mathcal {X}}$
is a triangulated category.

\begin{Prop}
   Let $A$ and $B$ be Gorenstein Artin algebras. Suppose that $F$
   is a derived equivalence between $A$ and $B$.
   Then we have the following statements.

   $(1)$ There is an equivalence
   $\underline{F}: \underline{_{A}\mathcal {X}}\ra\underline{_{B}\mathcal {X}}$.

   $(2)$ If $A$ and $B$ are finite dimensional algebras over a field
   $k$, then there exist bimodules $_{A}M_{B}$ and $_{B}L_{A}$ such
   that the pair of functors
   $$
_{A}M_{B}\otimes-: A\text{-}mod\ra B\text{-}mod,\;
_{B}L_{A}\otimes-: B\text{-}mod\ra A\text{-}mod
   $$
   induces an equivalence of triangulated categories
   $\underline{_{A}\mathcal {X}}$ and $\underline{_{B}\mathcal {X}}$.
\end{Prop}

\textbf{ Proof.} We refer to \cite[Theorem 4.6]{Ha2} and \cite
[Theorem 5.4]{Ka} for the proofs of (1) and (2), respectively.
$\square$

Our main result in this chapter is the following theorem.

\begin{Theo}\label{T} Let $A$ and $B$ be Gorenstein Artin algebras of
Cohen-Macaulay finite type. If $A$ and $B$ are derived equivalent,
then the Cohen-Macaulay Auslander algebras $\Lambda$ and $\Gamma$ of
$A$ and $B$ are also derived equivalent.
\end{Theo}

\textbf{ Proof.} In fact, if Artin algebras $A$ and $B$ are derived
equivalent, then $A$ is Gorenstein if and only if $B$ is Gorenstein.
By Proposition 3.10 or \cite[Theorem 8.11]{Be}, if Gorenstein Artin
algebras $A$ and $B$ are derived equivalent, then $A$ is of
Cohen-Macaulay finite type if and only if $B$ is. Let $F:\Db{A}\lra
\Db{B}$ be a derived equivalence. Set $\Lambda=\End(A\oplus X)$ with
$X=\oplus_{0\leq i\leq m} X_{i}$, where each $X_{i}$ is
indecomposable non-projective Gorenstein projective $A$-module. Then
$\Lambda$ is the Cohen-Macaulay Auslander algebra of $A$. By
Proposition 3.10, it follows that $Y_{i}=\underline{F}(X_{i})$ is
the indecomposable non-projective Gorenstein projective $B$-module.
Set $Y=\oplus_{0\leq i\leq m} Y_{i}$. Then $\Gamma=\End(B\oplus Y)$
is the Cohen-Macaulay Auslander algebra of $B$. Let $N$ be the
$B$-module $(B\oplus Y)$ and let $\cpx{\bar{T}}$ be the complex
$F(A\oplus X)$. Thus, we construct a tilting complex
$\Hom(N,\cpx{\bar{T}})$. The result follows from Proposition 3.8.
$\square$

\noindent{\bf Remark.} Let $A$ and $B$ be Gorenstein Artin algebras
of Cohen-Macaulay finite type. According to a result of Liu and Xi
\cite[Theorem 1.1]{LX2}, we see that, if $A$ and $B$ are stably
equivalent of Morita type, then the Cohen-Macaulay Auslander
algebras of $A$ and $B$ are also stably equivalent of Morita type.

As a corollary of Theorem \ref{T}, we re-obtain the following result
of Hu and Xi \cite{HX3} since self-injective Artin algebras of
finite representation type are Gorenstein Artin algebras of
Cohen-Macaulay finite type.

\begin{Koro} $\rm\cite[Corollary\,3.13]{HX3}$
Suppose that $A$ and $B$ are self-injective Artin algebras of finite
representation type. If $A$ and $B$ are derived equivalent, then the
Auslander algebras of $A$ and $B$ are also derived equivalent.
\end{Koro}

\noindent{\bf Acknowledgements.} The author would like to thank his
supervisor Professor Changchang Xi. He is grateful to him for his
guidance, patience and kindness.

\noindent{Shengyong Pan}\\

\noindent{Department of Mathematics,\\ Beijing Jiaotong University,
Beijing 100044,\\ People's Republic of China}\\

\noindent{School of Mathematical Sciences,\\
Beijing Normal University, Beijing 100875,\\ People's Republic of
China\\ \texttt{E-mail:shypan@bjtu.edu.cn}}

\bigskip
\today{}

\begin{thebibliography}{99}
{\small

\bibitem{AH}{{\sc H. Abe} and {\sc M. Hoshino},
Gorenstein orders associated with modules.\, {\it Comm. Algebra} 38
(2010), 165-180.}

\bibitem{ARS}{{\sc M. Auslander, I. Reiten} and {\sc S. O. Smal\o},
{\it Representation Thoery of Artin Algebras}. Cambridge University
Press, 1995.}

\bibitem{Be}{{\sc A. Beligiannis},
Cohen-Macaulay modules, (co)torsion pairs and virtually Goren- stein
algebras.\, {\it J. Algebra} 288(2005), 137-211. }


\bibitem{BR}{{\sc A. Beligiannis} and {\sc I. Reiten},
{\it Homological and Homotopical Aspects of Torsion Theories}. Mem.
Amer. Math. Soc., 188(2007).}

\bibitem{BMRR}{{\sc A. Buan, R. Marsh, M, Reineke} and {\sc I. Reiten},
Tilting theory and cluster combinatorics.\, Adv.\, Math. 204(2006),
572-618.}

\bibitem{Ch}{{\sc X. W. Chen},
Gorenstein homological algebra of Artin algebras, Postdoctoral
Report, USTC, 2010.}

\bibitem{FGR}{{\sc R. M. Fossum}, {\sc P. A. Griffith} and {\sc I. Riten},
{\it Trivial Extensions of Abelian Categories}.\, Springer Lecture
Notes 456, Heidelberg 1975.}

\bibitem{GZ}{{\sc N. Gao} and {\sc P. Zhang},
Gorenstein derived categories. {\it J. Algebra} 323(2010),
2041-2057}.

\bibitem{Ha1}{{\sc D. Happel},
{\it Triangulated Categories in the Representation Theory of Finite
Dimensional Algebras}. Cambridge University Press, Cambridge. 1988.}

\bibitem{Ha2}{{\sc D. Happel},
On Gorenstein algebras.\, {\it In: Representation theory of finite
groups and finitedimensional algebras }\, (Proc. Conf. at Bielefeld,
1991), 389-404, \, Progress in Math., vol. 95, Birkh$\ddot{a}$user,
Basel, 1991.}

\bibitem{HK1}{{\sc M. Hoshino} and {\sc Y. Kato},
Tilting complexes defined by idempotents. {\it Comm.\, Algebra}
30(2002), 83-100.}

\bibitem{HK2}{{\sc M. Hoshino} and {\sc Y. Kato},
Tilting complexes associated with a sequence of idempotents. {\it J.
Algebra} 183(2003), 105-124.}

\bibitem{Hu}{{\sc W. Hu},
On iterated almost $\nu$-stable derived equivalences. Preprint,
available at : http://arxiv.org/abs/0811.0926v3, 2008.}

\bibitem{HKX}{{\sc W. Hu}, {\sc S. Koenig} and {\sc C. C. Xi},
Derived eqivalences from cohomological approximations, and mutations
of perforated Yoneda algebras. Preprint, available at :
http://arxiv.org/abs/arXiv:1102.2790, 2011.}

\bibitem{HX1}{{\sc W. Hu} and {\sc C. C. Xi},
$\mathcal {D}$-split sequences and derived equivalences. To appear
in Adv. Math. (2011).}

\bibitem{HX2}{{\sc W. Hu} and {\sc C. C. Xi},
Derived equivalences and stable equivalences of Morita type, I.
Preprint, available at :
http://math.bnu.edu.cn/~ccxi/Papers/Articles/xihu-3.pdf, 2007, to
appear in Nagoya Math. J. }

\bibitem{HX3}{{\sc W. Hu} and {\sc C. C. Xi},
Derived equivalences for $\Phi$-Auslander-Yoneda algebras. Preprint,
available at :
http://math.bnu.edu.cn/~ccxi/Papers/Articles/xihu-4.pdf, 2009.}

\bibitem{Ka}{{\sc Y. Kato},
On derived equivalent coherent rings. \, {\it Comm.\, Algebra}
30(2002), 4437-4454.}

\bibitem{LX1}{{\sc Y. M. Liu} and {\sc C. C. Xi},
Constructions of stable equivalences of Morita type for
finite-dimensional algebras. II. \emph{Math. Z.} {\bf 251}(2005),
no.1, 21-39.}

\bibitem{LX2}{{\sc Y. M. Liu} and {\sc C. C. Xi},
Constructions of stable equivalences of Morita type for finite
dimensional algebras III. J. London Math. Soc. 76(2007), 567-585.}

\bibitem{N}{{\sc A. Neeman},
{\it Triangulated Categories}.\, Annals of Mathematics Studies 148,
Princeton University Press, \, Princeton and Oxford,\, 2001.}

\bibitem{P}{{\sc S. Y. Pan},
Derived equivalences for $\Phi$-Cohen Macaulay Auslander-Yoneda
algebras. In preparation.}


\bibitem{PX}{{\sc S. Y. Pan} and {\sc C. C. Xi},
Finiteness of finitistic dimension is invariant under derived
equivalences. \, {\it J. Algebra} 322(2009), 21-24.}

\bibitem{Ri1}{{\sc J. Rickard},
Morita theory for derived categories. {\it J. London Math. Soc.} 39
(1989), 436-456.}

\bibitem{Ri2}{{\sc J. Rickard},
Derived categories and stable equivalence. {\it J. Pure Appl.
Algebra} 61 (1989), 303-317.}

\bibitem{Ri3}{{\sc J. Rickard},
Derived equivalences as derived functors. {\it J. London Math. Soc.}
43 (1991), 37-48.}

\bibitem{Ver}{{\sc J. Verdier},
Cat\'{e}gories d\'{e}riv\'{e}es, \'{e}tat 0. \, {\it Lecture Notes
in Math. 569}(1977), Springer, Berlin, 262-311.} }
\end{thebibliography}
\end{document}